\documentstyle[twoside]{article}
\include{graphicx}
\title{\bf An asymptotic expansion for a sum of modified Bessel functions with quadratic argument}
\author{\sc R. B. Paris\footnote{E-mail address:\ \ {\tt r.paris@abertay.ac.uk}}\\
\\
{\em Division of Computing and Mathematics,}\\
{\em Abertay University, Dundee DD1 1HG, UK}\\
}

\begin{document}
\newcommand{\bee}{\begin{equation}}
\newcommand{\ee}{\end{equation}}
\def\f#1#2{\mbox{${\textstyle \frac{#1}{#2}}$}}
\def\dfrac#1#2{\displaystyle{\frac{#1}{#2}}}
\newcommand{\fr}{\frac{1}{2}}
\newcommand{\fs}{\f{1}{2}}
\newcommand{\g}{\Gamma}
\newcommand{\br}{\biggr}
\newcommand{\bl}{\biggl}
\newcommand{\ra}{\rightarrow}
\renewcommand{\topfraction}{0.9}
\renewcommand{\bottomfraction}{0.9}
\renewcommand{\textfraction}{0.05}
\newcommand{\mcol}{\multicolumn}
\date{}
\maketitle
\pagestyle{myheadings}
\markboth{\hfill {\it R.B. Paris} \hfill}
{\hfill {\it A modified Bessel function sum } \hfill}
\begin{abstract} 
We examine the sum of modified Bessel functions with argument depending quadratically on the summation index given by
\[S_\nu(a)=\sum_{n\geq 1} (\fs an^2)^{-\nu} K_\nu(an^2)\qquad (|\arg\,a|<\fs\pi)\]
as the parameter $|a|\to 0$. It is shown that the positive real $a$-axis is a Stokes line, where an infinite number of increasingly subdominant exponentially small terms present in the asymptotic expansion undergo a smooth, but rapid, transition as this ray is crossed. Particular attention is devoted to the details of the expansion on the Stokes line as $a\to 0$ through positive values. Numerical results are presented to support the asymptotic theory.
\vspace{0.4cm}

\vspace{0.3cm}

\end{abstract}

\vspace{0.2cm}

\noindent $\,$\hrulefill $\,$

\vspace{0.2cm}

\begin{center}
{\bf 1. \  Introduction}
\end{center}
\setcounter{section}{1}
\setcounter{equation}{0}
\renewcommand{\theequation}{\arabic{section}.\arabic{equation}}
We consider the asymptotic expansion of the sum
\bee\label{e11}
S_\nu(a)=\sum_{n\geq 1} (\fs an^2)^{-\nu} K_\nu(an^2)
\ee
as the parameter $a\to 0$ in $|\arg\,a|<\fs\pi$, where $K_\nu(z)$ is the modified Bessel function of the second kind
and the order $\nu\geq 0$. The argument of the Bessel function depends quadratically on the summation index $n$, thereby extending the sum considered in \cite{P17} with argument depending linearly on $n$. The sum converges without restriction on $\nu$ on account of the exponential decay of $K_\nu(an^2)\sim (\pi/2a)^{1/2}n^{-1} \exp [-an^2]$ as $n\to\infty$.

In the case $\nu=\fs$ we have $K_\frac{1}{2}(z)=e^{-z} (\pi/2z)^{1/2}$. The above sum then reduces to a case of the generalised Euler-Jacobi series $\sum_{n\geq 1}n^{-w} e^{ -an^p}$, $p>0$ and we find
\bee\label{e12}
S_\frac{1}{2}(a)=\frac{\sqrt{\pi}}{a} \sum_{n\geq 1}\frac{e^{-an^2}}{n^2}.
\ee
Sums of this type were first studied by Ramanujan  who showed that the expansion as $a\to 0$ consisted of an asymptotic sum involving the Riemann zeta function; see \cite[Chapter 1]{Berndt}. A detailed study using a hypergeometric function approach for the Euler-Jacobi series when $w=0$ and $p$ is a rational fraction has been given in the monograph \cite{Kow}, where great emphasis was placed on obtaining exponentially small contributions. More  recently, the generalised Euler-Jacobi series has been considered in \cite{P16}; see also \cite[Section 8.1]{PK} for the case with $w=0$ and $p>0$. 

The sum $S_\nu(a)$ becomes difficult to compute in the limit $a\to 0$ due to the resulting slow convergence of the series. It will be shown that the asymptotic expansion of $S_\nu(a)$ consists of a finite series in ascending powers of $a$ with coefficients involving the Riemann zeta function, together with an infinite number of increasingly subdominant exponentially small terms. These exponential terms experience a Stokes phenomenon as the positive real $a$-axis is crossed. Each exponential is associated with a multiplier that undergoes a smooth, but rapid, change described to leading order by an error function of appropriate argument from the value $\exp [-\pi i(\nu-\fs)]$ just above the axis to $\exp [\pi i(\nu-\fs)]$ just below the axis. Other well-known functions that exhibit an infinite number of subdominant exponential terms in their expansions as $|z|\to\infty$ 
are $\log\,\g(z)$ \cite[\S 6.4]{PK} and the Hurwitz zeta function $\zeta(s,z)$ \cite{P05}.

We pay particular attention to the expansion of $S_\nu(a)$ as $a\to 0$ on the positive real axis. In this case the multipliers of the exponentially small terms are not given by $\sin \pi\nu$, as might appear from the above, but receive additional contributions from other exponentially small terms that become comparable on the positive real $a$-axis.

\vspace{0.6cm}

\begin{center}
{\bf 2. \ An expansion for $S_\nu(a)$}
\end{center}
\setcounter{section}{2}
\setcounter{equation}{0}
\renewcommand{\theequation}{\arabic{section}.\arabic{equation}}
We start with the Mellin-Barnes integral representation\footnote{There is an error in the sector of validity in \cite[(10.32.13)]{DLMF} and also in \cite[p.~114]{PK}.} \cite[(10.32.13)]{DLMF}
\[(\fs x)^{-\nu} K_\nu(x)=\frac{1}{4\pi i}\int_{c-\infty i}^{c+\infty i}\g(s) \g(s-\nu) (\fs x)^{-2s} ds\qquad (|\arg\,x|<\fs\pi),\]
where $c>\max\{0,\nu\}$ so that the integration path lies to the right of all the poles of the integrand. The quantity $c$ will be used as a generic parameter that can vary according to each integral. Then it follows that
\[
S_\nu(a)=\sum_{n\geq 1}(\fs an^2)^{-\nu} K_\nu(an^2)=\frac{1}{4\pi i}\int_{c-\infty i}^{c+\infty i} \g(s)\g(s-\nu)(\fs a)^{-2s} \sum_{n\geq 1} n^{-4s}\,ds\]
\bee\label{e21}
=\frac{1}{8\pi i}\int_{c-\infty i}^{c+\infty i} \g(\fs s) \g(\fs s-\nu) (\fs a)^{-s} \zeta(2s)\,ds\qquad (|\arg\,a|<\fs\pi),
\ee
where $\zeta(s)$ is the Riemann zeta function and $c>\max\{\fs, 2\nu\}$. Throughout we let $\theta=\arg\,a$.

Let $\nu=m+\nu'$, where $m=0, 1, 2, \ldots$ and $0\leq\nu'<1$. The integrand in (\ref{e21}) has a pole at $s=\fs$ resulting from $\zeta(2s)$; it also has poles at $s=2\nu-2n$, $n=0, 1, 2, \ldots$ together with a pole at $s=0$, the other poles of $\g(\fs s)$ at $s=-2, -4, \ldots$ being cancelled by the trivial zeros of $\zeta(2s)$ at $s=-1, -2, \ldots\, $. It may be noted that when $2\nu$ is an odd integer (that is, when $\nu'=\fs$) the sequence of poles at $s=2\nu-2n$ is finite with $0\leq n\leq m$. This case is related to the Euler-Jacobi series studied in \cite{P16} since, when $\nu=m+\fs$, $K_\nu(x)$ can be expressed in terms of exponentials.

We consider the integral taken round the rectangular contour with vertices at $c\pm iT$, $-d\pm iT$, where $0<d<2(1-\nu')$.
The contribution from the upper and lower sides $s=\sigma\pm iT$, $-d\leq \sigma\leq c$ as $T\to\infty$ can be estimated by use of the standard results
\[|\g(\sigma\pm it)|\sim \sqrt{2\pi}t^{\sigma-\frac{1}{2}}e^{-\frac{1}{2}\pi t}\qquad (t\to+\infty),\]
which follows from Stirling's formula for the gamma function, and \cite[p.~25]{Iv}
\[|\zeta(\sigma\pm it)|=O(t^{{\hat\mu}(\sigma)}\log^\beta t)\qquad (t\to+\infty)\]
where ${\hat\mu}(\sigma)=0$ ($\sigma>1$), $\fs-\fs\sigma$ ($0\leq\sigma\leq1$), $\fs-\sigma$ ($\sigma\leq0$) and $\beta=0$ ($\sigma>1$), 1 ($\sigma\leq 1)$. Then it follows that
\[\g(\fs\sigma\pm\fs it) \g(\fs\sigma-\nu\pm\fs it)(\fs a)^{-\sigma\mp it}=O(t^{\sigma-\nu-1} e^{-\frac{1}{2}\pi t+|\theta|t})\]
as $t\to +\infty$.
Hence the modulus of the integrand on these horizontal paths is $O(T^{\xi} \log\,T e^{-\Delta T})$ as $T\to\infty$, where $\xi=\sigma+{\hat\mu}(\sigma)-\nu-1$, $\Delta=\fs\pi-|\theta|$. Taking into account the different forms of ${\hat\mu}(\sigma)$, we see that the contribution from these paths vanishes as $T\to\infty$ provided $|\theta|<\fs\pi$.

We first displace the integration path to the left over the poles in $\Re (s)\geq 0$. Three situations can arise: the poles at $s=0$, $\fs$ and $2\nu-2n$, $0\leq n\leq m$ can (i) be all simple ($\nu'\neq 0, \f{1}{4}$), (ii) be all simple save for a double pole at $s=0$ ($\nu'=0$) and (iii) be all simple save for a double pole at $s=\fs$ ($\nu'=\f{1}{4}$). Making use of the results that $\zeta(s)$ has residue equal to 1 at $s=1$, $\zeta(0)=-\fs$, $\zeta'(0)=-\fs\log\,2\pi$ and
\[\g(\alpha+\epsilon)=\g(\alpha)\{1+\epsilon \psi(\alpha)+O(\epsilon^2)\},\qquad \zeta(1+\epsilon)=\epsilon^{-1}\{1+\epsilon\gamma+O(\epsilon^2)\} \qquad (\epsilon\to 0),\]
where $\psi(\alpha)$ denotes the logarithmic derivative of the gamma function and $\gamma=0.577215\ldots$ is the Euler-Mascheroni constant, we then find after routine calculations the residue contribution in $\Re (s)\geq 0$ given by
\bee\label{e2h}
{\cal H}(a;\nu)=\left\{\begin{array}{ll} \displaystyle{\frac{1}{2}\sum_{n=0}^m h_n\!+\!\frac{\g(\f{1}{4})\g(\f{1}{4}-\nu)}{8\surd(\fs a)} \!+\!\frac{\pi\mbox{cosec} \pi\nu}{4\g(1+\nu)}} & (\nu'\neq 0, \f{1}{4})\\
\\
\displaystyle{\frac{1}{2}\sum_{n=0}^{m-1} h_n\!+\!\frac{\g(\f{1}{4})\g(\f{1}{4}-\nu)}{8\surd(\fs a)}\!+\!\frac{(-)^m}{2 m!}\{\fs\gamma-\fs\psi(m+1)+\log \frac{a}{8\pi^2}\}} & (\nu'=0)\\
\\
\displaystyle{\frac{1}{2}\sum_{n=0}^{m-1} h_n \!+\!\frac{\pi\mbox{cosec} \pi\nu}{4\g(1+\nu)}\!+\!\frac{(-)^m\g(\f{1}{4})}{4m! \surd(\fs a)}  \{\fs\psi(\f{1}{4})\!+\!\fs\psi(m+1)\!+\!2\gamma\!-\!\log \fs a\}}
& (\nu'=\f{1}{4}),
\end{array}\right.\ee
where
\[h_n:=\frac{(-)^n }{n!} \g(\nu-n)\zeta(4\nu-4n) (\fs a)^{2n-2\nu}.\]
Then we have
\bee\label{e23}
S_\nu(a)-{\cal H}(a;\nu)=\frac{1}{8\pi i}\int_{-c-\infty i}^{-c+\infty i}\g(\fs s) \g(\fs s-\nu) (\fs a)^{-s} \zeta(2s)\,ds,
\ee
where $0<c<2(1-\nu')$. 

We now employ the functional relation for $\zeta(s)$ in the form \cite[(25.4.2)]{DLMF}
\[\zeta(s)=2^s\pi^{s-1} \g(1-s)\zeta(1-s) \sin \fs\pi s\]
and make the change of variable $s\to -s$ to obtain
\[S_\nu(a)-{\cal H}(a;\nu)=-\frac{\pi}{8\pi i}\int_{c-\infty i}^{c+\infty i} \frac{\g(1+2s)\zeta(1+2s)}{\g(1\!+\!\fs s)\g(1\!+\!\fs s\!+\!\nu)}\,\bl(\frac{a}{8\pi^2}\br)^{\!\!s}\,\frac{\sin \pi s}{\sin \fs\pi s\,\sin \pi(\fs s\!+\!\nu)}\,ds\]
\[=\sum_{k\geq 1}\frac{1}{k}\,J_k(a;\nu),\]
where
\bee\label{e24}
J_k(a;\nu)=-\frac{\pi}{4\pi i}\int_{c-\infty i}^{c+\infty i}
\frac{\g(1\!+\!2s)}{\g(1\!+\!\fs s)\g(1\!+\!\fs s\!+\!\nu)}\,
\bl(\frac{a}{8\pi^2k^2}\br)^{\!\!s}\,\frac{\cos \fs\pi s}{\sin \pi(\fs s+\nu)}\,ds
\ee
with $0<c<2(1-\nu')$.
Here we have used the series expansion for $\zeta(1+2s)$ since $\Re (s)>0$ on the integration path
and reversal of the order of summation and integration is justified when $|\arg\,a|<\fs\pi$ by absolute convergence.

Displacement of the integration path in (\ref{e24}) to the right over the simple poles at $s=2n-2\nu$, $m+1\leq n\leq N_k-1$, where $N_k$ is (for the moment) an arbitrary positive integer, then yields
\bee\label{e25}
J_k(a;\nu)=\cos \pi\nu \!\!\!\sum_{n=m+1}^{N_k-1} \frac{(4n-4\nu)!}{n! (n-\nu)!}\,\bl(\frac{a}{8\pi^2 k^2}\br)^{\!2n-2\nu}
+R_k(a;N_k).
\ee
With the further change of variable $s\to s-2\nu$, the remainder is given by
\bee\label{e26}
R_k(a;N_k)=-\frac{\pi}{4\pi i}\int_{{\cal L}_N}
\frac{\g(1\!+\!2s\!-\!4\nu)}{\g(1\!+\!\fs s) \g(1\!+\!\fs s\!-\!\nu)}\,\bl(\frac{a}{8\pi^2 k^2}\br)^{\!s-2\nu}\, \Theta(s)\,ds,
\ee
where
the integration path ${\cal L}_N$ denotes the rectilinear path $(-c\!+\!2N_k\!-\!\infty i, -c\!+\!2N_k\!+\!\infty i)$, with $0<c<2$, and
\[\Theta(s):=\frac{\cos \fs\pi(s-2\nu)}{\sin \fs\pi s}=\frac{\cos \pi\nu(1+\cos \pi s)}{\sin \pi s}+\sin \pi\nu.\]

Then we have the result
\bee\label{e27}
S_\nu(a)={\cal H}(a;\nu)+\sum_{k\geq 1}\frac{1}{k}\bl\{\cos \pi\nu\!\!\!\sum_{n=m+1}^{N_k-1} \frac{(4n-4\nu)!}{n! (n-\nu)!}\,\bl(\frac{a}{8\pi^2 k^2}\br)^{\!2n-2\nu}\!\!
+R_k(a;N_k)\br\}.
\ee
The contribution to $S_\nu(a)$ from $\Re (s)<0$ in the integral in (\ref{e23}) has been decomposed into a $k$-sequence of component asymptotic series with scale $8\pi^2k^2/a$, each associated with its own truncation index $N_k$ and remainder term $R_k(a;N_k)$. The indices $N_k$ ($k=1, 2, \ldots$) are, for the moment, unspecified but will be suitably chosen in the next section.
\vspace{0.6cm}

\begin{center}
{\bf 3. \ Expansion of the remainder term $R_k(a;N_k)$}
\end{center}
\setcounter{section}{3}
\setcounter{equation}{0}
\renewcommand{\theequation}{\arabic{section}.\arabic{equation}}
We now choose the indices $N_k$ to correspond to the optimal truncation value of the $k$-component asymptotic series in (\ref{e25}) (that is, truncation just before the numerically smallest term). This is given by (when we assume $N_k\gg 1$)
\[\bl(\frac{8\pi^2k^2}{|a|}\br)^{\!2} \simeq \frac{(4N_k+4-4\nu)!}{(4N_k-4\nu)!}\,\frac{1}{(N_k+1)(N_k+1-\nu)}\] 
\[=\frac{(4N_k)^4}{N_k^2}\,\frac{\prod_{r=1}^4(1-\frac{r-4\nu}{4N_k})}{(1+\frac{1}{N_k})(1+\frac{1-\nu}{N_k})}=(16N_k)^2\bl\{1+\frac{1-6\nu}{2N_k}+O(N_k^{-2})\br\}.\]
If we define the quantities
\bee\label{e31}
X_k:=\frac{\pi^2k^2}{a},\qquad \vartheta:=-3\nu,
\ee
then the optimal truncation indices are specified by
\bee\label{e32}
|X_k|=2N_k+\vartheta+\alpha,
\ee
where $|\alpha|$ is bounded. It is seen that as $|a|\to 0$ the truncation indices $N_k\to+\infty$, ($k=1,2, \ldots$).
 
Since $N_k\to+\infty$, the variable $s$ in the quotient of gamma functions in (\ref{e26}) is uniformly large on the displaced integration path ${\cal L}_N$. From Lemma 2.2 in \cite[p.~39]{PK}, we have the inverse factorial expansion
\bee\label{e3e}
\frac{\g(1+2s-4\nu)}{\g(1\!+\!\fs s) \g(1\!+\!\fs s\!-\!\nu)}=\frac{2^{3s+\frac{3}{2}-5\nu}}{2\pi}\bl\{\sum_{j=0}^{M-1}(-)^j c_j(\nu) \g(s+\vartheta-j)+\rho_M(s) \g(s+\vartheta-M)\br\}
\ee
for positive integer $M$. The remainder function $\rho_M(s)$ is analytic in $s$ except at the points $s=(4\nu-1-r)/2$, $r=0, 1, 2 \ldots$ and is such that $\rho_M(s)=O(1)$ as $|s|\to\infty$ in $|\arg\,s|<\pi$. The coefficients $c_j(\nu)$ are discussed in the appendix where their values for $j\leq 4$ are presented. We have the values for $j\leq 3$
\[c_0(\nu)=1,\quad c_1(\nu)=\frac{3}{8}(1+2\nu)^2,\]
\bee\label{e33a}
c_2(\nu)=\frac{1}{128}(57+344\nu+696\nu^2+544\nu^3+144\nu^4),
\ee
\[c_3(\nu)=\frac{1}{1024}(945+7068\nu+19196\nu^2+24288\nu^3+15472\nu^4+4800\nu^5+576\nu^6).\] 
Then, we obtain
\bee\label{e33}
R_k(a;N_k)=-2^{\nu-\frac{1}{2}}\,\bl(\frac{a}{\pi^2k^2}\br)^{\!\nu}\bl\{\sum_{j=0}^{M-1}\frac{(-)^j c_j(\nu)}{X_k^j}\,\frac{1}{2\pi i}\int_{{\cal L}_N} \g(s\!+\!\vartheta\!-\!j) X_k^{-s-\vartheta+j}\Theta(s)\,ds+{\cal R}_{M,N_k}\br\},
\ee
where
\[{\cal R}_{M,N_k}=\frac{1}{2\pi i}\int_{{\cal L}_N}\rho_M(s) \g(s\!+\!\vartheta\!-\!M) X_k^{-s-\vartheta} \Theta(s)\,ds.\] 
 
The remainder ${\cal R}_{M,N_k}$ can be estimated by application of the bounds discussed in \cite[Section 2.5.2]{PK}.
With $\Theta(s)=\Theta_1(s)+\sin \pi\nu$, the part of the remainder integral containing $\Theta_1(s)$ can be
split into three separate integrals with variables $X_k^{-s}$ and $(X_ke^{\pm \pi i})^{-s}$, to each of which we can apply Lemma 2.9 in \cite[p.~75]{PK}. We therefore obtain the order estimates as $|X_k|\to\infty$ in  $|\arg\,X_k|<\fs\pi$ given by $O(X_k^{-M-1/2}e^{-|X_k|})$ for the integral with variable $X_k^{-s}$ and $O(X_k^{-M}e^{-|X_k|})$ for the integrals\footnote{In \cite[(2.5.11)]{PK} the bound is obtained only in the sector $|\arg\,z|\leq\pi$. The domain of validity is easily extended to $|\arg\,z|<\frac{3}{2}\pi$ by application of Lemma 2.6 in \cite[p.~70]{PK}.} with variables $(X_ke^{\pm\pi i})^{-s}$. For the part of the remainder integral containing $\sin \pi\nu$, we have from Lemma 2.7 in \cite[p.~71]{PK} the estimate $O(X_k^{-M}e^{-X_k})$.
Hence
${\cal R}_{M,N_k}=O(X_k^{-M}e^{-X_k})$ as $|X_k|\to+\infty$ in $|\arg\,X_k|<\fs\pi$.
 
We now introduce the so-called {\it terminant function} $T_\nu(z)$ defined\footnote{In \cite[(2.11.11)]{DLMF} this function is denoted by $F_\nu(z)$ and is expressed as a multiple of the exponential integral $E_\nu(z)=z^{\nu-1}\g(1-\nu,z)$. In addition, the parameter $\nu$ here is not to be confused with that appearing in (\ref{e11}).} as a multiple of the incomplete gamma function $\g(a,z)$ by
\[{T}_\nu(z):=\frac{\g(\nu)}{2\pi}\, \g(1-\nu,z).\]
From the formula connecting $\g(a,ze^{\pm\pi i})$ given in \cite[(8.2.10)]{DLMF}
we have the connection formula (compare also \cite[(6.2.45)]{PK})
\bee\label{e34}
T_\nu(ze^{-\pi i})=e^{2\pi i\nu}\{T_\nu(ze^{\pi i})-ie^{-\pi i\nu}\}.
\ee
The Mellin-Barnes integral representation of this function is \cite[(6.2.7)]{PK}
\bee\label{e35}
-2 e^z T_\nu(z)=\frac{1}{2\pi i}\int_{-c-\infty i}^{-c+\infty i} \frac{\g(s+\nu)}{\sin \pi s}\,z^{-s-\nu}ds\qquad (|\arg\,z|<\f{3}{2}\pi,\ 0<c<1)
\ee
provided $\nu\neq 0, -1, -2, \ldots\,$. 

Then, if we make the change of variable $s\to s+2N_k$ in the integrals appearing in (\ref{e33}), we obtain 
\bee\label{e36}
R_k(a;N_k)=2^{\nu-\frac{1}{2}} \bl(\frac{a}{\pi^2k^2}\br)^{\!\nu} \bl\{\sum_{j=0}^{M-1}\frac{(-)^j c_j(\nu)}{X_k^j}\,\Upsilon_j(a)+O(X_k^{-M}e^{-X_k})\br\},
\ee
where 
\[\Upsilon_j(a)=-\frac{1}{2\pi i}\int_{-c-\infty i}^{-c+\infty i} \g(s\!+\!\mu\!-\!j) X_k^{-s-\mu+j}\bl\{\frac{\cos \pi\nu (1\!+\!\cos \pi s)}{\sin \pi s}+\sin \pi\nu\br\}ds, \quad \mu:=2N_k+\vartheta.\]
This last integral can be evaluated using (\ref{e35}), together with the Cahen-Mellin integral \cite[p.~90]{PK} for the part involving $\sin \pi\nu$, to yield upon writing $\cos \pi s$ in terms of exponentials 
\bee\label{e37}
\Upsilon_j(a)=e^{-X_k}\{2\cos \pi\nu\,{\bf T}(X_k)-\sin \pi\nu\},
\ee
where
\[{\bf T}(X_k)=e^{2X_k} T_{\mu-j}(X_k)
+\frac{(-)^j}{2}\,\{e^{\pi i\vartheta} T_{\mu-j}(X_k e^{\pi i})+e^{-\pi i\vartheta} T_{\mu-j}(X_k e^{-\pi i})\}\]
\[=e^{2X_k} T_{\mu-j}(X_k)+(-)^j e^{\pi i\vartheta} T_{\mu-j}(X_k e^{\pi i})-\fs i\]
upon use of the connection formula (\ref{e34}) in the form
\[e^{-\pi i\vartheta} T_{\mu-j}(X_ke^{-\pi i})=e^{\pi i\vartheta} T_{\mu-j}(X_ke^{\pi i})-(-)^ji.\]

Combination of (\ref{e36}) and (\ref{e37}) then provides an expression for $R_k(a;N_k)$ in terms of terminant functions in the form
\bee\label{e38}
R_k(a;N_k)=2^{\nu-\frac{1}{2}} \bl(\frac{a}{\pi^2 k^2}\br)^{\!\nu} e^{-X_k} \bl\{\sum_{j=0}^{M-1} \frac{(-)^j c_j(\nu)}{X_k^j}\bl(2\cos \pi\nu \,{\bf T}(X_k)-\sin \pi\nu\br)
+O(X_k^{-M})\br\},
\ee
It now remains to exploit the known asymptotics of these terminant functions when $\mu \sim |X_k|$ as $|X_k|\to\infty$.

\vspace{0.6cm}

\begin{center}
{\bf 4. \  The asymptotic expansion for $R_k(a;N_k)$}
\end{center}
\setcounter{section}{4}
\setcounter{equation}{0}
\renewcommand{\theequation}{\arabic{section}.\arabic{equation}}
The asymptotic expansion of the terminant function $T_\nu(z)$ for large $\nu$ and complex $z$, when $\nu\sim |z|$, has been discussed in detail by Olver in \cite{O91}; see also \cite[Section 2.11]{DLMF} and the detailed account in \cite[pp.~259--265]{PK}. It is found that
\[T_\nu(z)=\left\{\begin{array}{ll}
\displaystyle{\frac{e^{-\phi\nu i}}{1+e^{-i\phi}}\,\frac{e^{-z-|z|}}{\sqrt{2\pi |z|}}\{1+O(z^{-1})\}} & |\arg\,z|\leq \pi-\delta, \delta>0\\
\displaystyle{ie^{-\pi i\nu}\{\fs+\fs \mbox{erf}\,[c(\phi) (\fs|z|)^{1/2}]\}+O\bl(\frac{e^{-z-|z|}}{\sqrt{2\pi|z|}}\br)} & \delta\leq\arg\,z\leq2\pi-\delta \end{array}\right.\]
valid as $|z|\to\infty$ when $\nu\sim |z|$, where $\phi=\arg\,z$ and the quantity $c(\phi)=\phi-\pi+O((\phi-\pi)^2)$ in the neighbourhood of $\phi=\pi$. 

Then with $\theta=\arg\,a$, $X_k=|X_k| e^{-i\theta}$ we have from (\ref{e37}) (with $\phi=\pi-\theta$)
\[\Upsilon_j(a)=-e^{-X_k} \sin \pi\nu+2\cos \pi\nu \bl\{O\bl(\frac{e^{-|X_k|}}{\sqrt{2\pi|X_k|}}\br)+e^{-X_k}\bl[\frac{1}{2}i\,\mbox{erf}\,[c(\phi) (\fs|X_k|)^{1/2}]\]
\[\hspace{4cm}+O\bl(\frac{e^{-X_k-|X_k|}}{\sqrt{2\pi|X_k|}}\br)\br]\br\}\]
\bee\label{e41}
=e^{-X_k} \bl(i \cos \pi\nu \ \mbox{erf}\,[c(\phi) (\fs|X_k|)^{1/2}]-\sin \pi\nu\br)+O\bl(\frac{e^{-|X_k|}}{\sqrt{2\pi|X_k|}}\br).
\ee
The form (\ref{e41}) approximately describes the smooth change of the multiplier associated with each exponential $e^{-X_k}$ --- on the increasingly sharp scale $(\fs|X_k|)^{1/2}=\pi k/\sqrt{2a}$ --- in the vicinity of the positive real $a$-axis. 
Noting that $c(\phi)\simeq -\theta$ near $\theta=0$, we see that the multiplier changes smoothly from approximately $-\sin \pi\nu-i\cos \pi\nu=e^{-\pi i(\nu+\frac{1}{2})}$
when $\theta>0$ to approximately $-\sin \pi\nu+i\cos \pi\nu=e^{\pi i(\nu+\frac{1}{2})}$ when $\theta<0$. 

The value of the multiplier when $\theta=0$ is not equal to $-\sin \pi\nu$, as would be expected from the above considerations, since additional exponentially small terms come into play on the positive real $a$-axis. We consider the expansion of $R_k(a;N_k)$ as $a\to 0$ through positive values in the next sub-section.
\vspace{0.3cm}

\noindent {\bf 4.1\ \ The expansion of $R_k(a;N_k)$ for $a>0$.}\ \ 
To deal with the expansion of $R_k(a;N_k)$, and hence $S_\nu(a)$, when $a\to 0$ with $\arg\,a=0$, we require more precise asymptotic information on the terminant function than that given above.
By expressing $T_\nu(z)$ in terms of the Laplace integral
\[T_\nu(z)=\frac{e^{-z}}{2\pi}\int_0^\infty e^{-zt}\,\frac{t^{\nu-1}}{1+t}\, dt,\]
Olver \cite{O91} established by application of the saddle-point method that 
\bee\label{e42}
T_{\mu-j}(x)=\frac{e^{-2x}}{2\sqrt{2\pi x}}\bl\{\sum_{r=0}^{K-1}A_{r,j}x^{-r}+O(x^{-K})\br\}\qquad (x\to+\infty),
\ee
when $\mu\sim x$ (and bounded integer $j$),
where we recall that $\mu=2N_k+\vartheta$. The coefficients $A_{r,j}$ satisfy $A_{0,j}=1$ ($j\geq 0$) and
\begin{eqnarray}
A_{1,j}&=&\f{1}{6}(2-6\gamma_j+3\gamma_j^2),\quad A_{2,j}=\f{1}{288}(-11-120\gamma_j+300\gamma_j^2-192\gamma_j^3+36\gamma_j^4),\nonumber\\
A_{3,j}&=&\f{2}{51840}(-587+3510\gamma_j+9765\gamma_j^2-26280\gamma_j^3+18900\gamma_j^4-5400\gamma_j^5+540\gamma_j^6),\nonumber\\
A_{4,j}&=&\f{1}{2448320}(120341-44592\gamma_j-521736\gamma_j^2-722880\gamma_j^3+2336040\gamma_j^4-1826496\gamma_j^5\nonumber\\
&&\hspace{5cm}+635040\gamma_j^6-103680\gamma_j^7+6480\gamma_j^8),\label{e31a}
\end{eqnarray}
where 
\bee\label{e42b}
\gamma_j:=\mu-x-j\qquad (0\leq j\leq K-1).
\ee

On the negative real $z$-axis, where a saddle point and a simple pole become coincident in the above Laplace integral, we have the expansion \cite{O91}
\bee\label{e43}
(-)^j e^{\pi i\vartheta}T_{\mu-j}(xe^{\pi i})=\frac{1}{2} i+\frac{1}{\sqrt{2\pi x}}\bl(\sum_{r=0}^{K-1} (\fs)_r G_{2r,j}\,(\fs x)^{-r}+O(x^{-K})\br)\qquad (x\to+\infty),
\ee
where the coefficients $G_{r,j}$ result from the expansion
\[\frac{\tau^{\gamma_j-1}}{1-\tau}\,\frac{d\tau}{dw}=-\frac{1}{w}+\sum_{r=0}^\infty G_{r,j}w^r,\qquad \fs w^2=\tau-\log\,\tau-1.\]
The branch of $w(\tau)$ is chosen such that $w\sim \tau-1$ as $\tau\ra 1$. Upon reversion of the $w$-$\tau$ mapping to yield
\[\tau=1+w+\f{1}{3}w^2+\f{1}{36}w^3-\f{1}{270}w^4+\f{1}{4320}w^5+ \cdots\ ,\]
it is found with the help of {\it Mathematica} that the first five even-order coefficients $G_{2r,j}\equiv 6^{-2r} {\hat G}_{2r,j}$ are
\begin{eqnarray}
{\hat G}_{0,j}\!\!&=&\!\!\f{2}{3}-\gamma_j,\qquad {\hat G}_{2,j}=\f{1}{15}(46-225\gamma_j+270\gamma_j^2-90\gamma_j^3), \nonumber\\
{\hat G}_{4,j}\!\!&=&\!\!\f{1}{70}(230-3969\gamma_j+11340\gamma_j^2-11760\gamma_j^3+5040\gamma_j^4
-756\gamma_j^5),\nonumber\\
{\hat G}_{6,j}\!\!&=&\!\!\f{1}{350}(-3626-17781\gamma_j+183330\gamma_j^2-397530\gamma_j^3+370440\gamma_j^4
-170100\gamma_j^5\nonumber\\
&&\hspace{7cm}+37800\gamma_j^6-3240\gamma_j^7),\nonumber\\
{\hat G}_{8,j}\!\!&=&\!\!\f{1}{231000}(-4032746+43924815\gamma_j+88280280\gamma_j^2-743046480\gamma_j^3\nonumber\\
&&+1353607200\gamma_j^4-1160830440\gamma_j^5+541870560\gamma_j^6
-141134400\gamma_j^7\nonumber\\
&&\hspace{6cm}+19245600\gamma_j^8-1069200\gamma_j^9).\label{e44}
\end{eqnarray}

Substitution of the expansions (\ref{e42}) and (\ref{e43}) into (\ref{e37}), with $x$ replaced by $X_k$ and $\gamma_j=\mu-X_k-j=-\alpha-j$ by (\ref{e32}),  then yields
\[\Upsilon_j(a)=e^{-X_k}\bl\{\frac{\cos \pi\nu}{\sqrt{2\pi X_k}}\bl( \sum_{r=0}^{K-1} D_{r,j} X_k^{-r}+O(X_k^{-K})\br)-\sin \pi\nu\br\}\]
as $a\to 0+$,
where
\bee\label{e45a}
D_{r,j}\equiv D_{r,j}(\alpha)=A_{r,j}+2^{r+1}(\fs)_r G_{2r,j}.
\ee
From (\ref{e36}) we therefore find (where we put $K=M$ for convenience)
\[R_k(a;N_k)=2^{\nu-\frac{1}{2}} \bl(\frac{a}{\pi^2 k^2}\br)^{\!\nu} e^{-X_k}\bl\{\sum_{j=0}^{M-1}\frac{(-)^j c_j(\nu)}{X_k^j}\bl(\frac{\cos \pi\nu}{\sqrt{2\pi X_k}} \sum_{r=0}^{M-1} D_{r,j} X_k^{-r}-\sin \pi\nu\br)+O(X_k^{-M})\br\}\]
\[=2^{\nu-\frac{1}{2}} \bl(\frac{a}{\pi^2 k^2}\br)^{\!\nu} e^{-X_k}\bl\{\frac{\cos \pi\nu}{\sqrt{2\pi X_k}}\sum_{j=0}^{M-1}(-)^j B_j X_k^{-j}-\sin \pi\nu \sum_{j=0}^{M-1} (-)^j c_j(\nu) X_k^{-j}+O(X_k^{-M})\br\},\]
where the coefficients $B_j$ are given by
\bee\label{e45}
B_j\equiv B_j(\alpha,\nu)=\sum_{r=0}^j (-)^r c_{j-r}(\nu) D_{r,j-r}.
\ee

Then we obtain the following theorem:
\newtheorem{theorem}{Theorem}
\begin{theorem}$\!\!\!.$\  Let $\nu=m+\nu'$, with $m=0, 1, 2, \ldots\ $ and $0\leq\nu'<1$. Further, let $X_k=\pi^2 k^2/a$ and the optimal truncation indices $N_k$ $(k=1, 2, \ldots)$ satisfy $X_k=2N_k+\vartheta+\alpha$, where $\vartheta=-3\nu$ and $|\alpha|$ is bounded. Then we have the expansion
\[S_\nu(a)={\cal H}(a;\nu)+\cos \pi\nu \sum_{k\geq 1}\frac{1}{k}\sum_{n=m+1}^{N_k-1} \frac{(4n-4\nu)!}{n! (n-\nu)!}\,\bl(\frac{a}{8\pi^2 k^2}\br)^{\!2n-2\nu}\]
\bee\label{e46}
+\frac{2^{\nu-\frac{1}{2}}a^\nu}{\pi^{2\nu}}\sum_{k\geq 1} \frac{e^{-X_k}}{k^{1+2\nu}}\bl\{\frac{\cos \pi\nu}{\sqrt{2\pi X_k}} \sum_{j=0}^{M-1}(-)^j B_j X_k^{-j}-\sin \pi\nu \sum_{j=0}^{M-1} (-)^j c_j(\nu) X_k^{-j}+O(X_k^{-M})\br\}
\ee
as $a\to 0+$. The sum ${\cal H}(a;\nu)$ is defined in (\ref{e2h}) and the coefficients $c_j(\nu)$ and $B_j$ are defined in (\ref{e33a}) and (\ref{e45}).
\end{theorem}

In Appendix B it is shown that the first double sum on the right-hand side of (\ref{e46}) can be rearranged in the form
\bee\label{e46a}
\cos \pi\nu \sum_{k\geq1} \sum_{n=N_{n-1}}^{N_n-1} \frac{(4n-4\nu)!}{n! (n-\nu)!} \bl(\frac{a}{8\pi^2}\br)^{\!2n-2\nu}\!\zeta(4n\!-\!4\nu\!+\!1,k),
\ee
where $N_0=m+1$ and $\zeta(x,k)$ is the Hurwitz zeta function \cite[p.~607]{DLMF}.
\vspace{0.3cm}

\noindent {\bf 4.2\ \ Two special cases}\ \ When $\nu=0$ we have from (\ref{e2h})
\[{\cal H}(a;0)=\frac{\g^2(\f{1}{4})}{2^{5/2} a^{1/2}}+\frac{1}{2}\bl\{\gamma+\log \frac{a}{8\pi^2}\br\}.\]
Then from (\ref{e46}) we obtain the expansion
\[S_0(a)={\cal H}(a;0)+\sum_{k\geq 1} \frac{1}{k} \sum_{n=1}^{N_k-1} \frac{(4n)!}{(n!)^2}\bl(\frac{a}{8\pi^2 k^2}\br)^{\!2n}
\!\!+\frac{a^{1/2}}{2\pi^{3/2}}\sum_{k\geq 1}\frac{e^{-X_k}}{k^2}\bl\{\sum_{j=0}^{M-1}(-)^j B_jX_k^{-j}+O(X_k^{-M})\br\}\]
as $a\to 0+$.

When $\nu=\fs$ we have
\[{\cal H}(a;\fs)=\frac{\pi^{5/2}}{6a}-\frac{\pi}{a^{1/2}}+\frac{1}{2}\pi^{1/2}\]
and 
\[S_\frac{1}{2}(a)={\cal H}(a;\fs)-\frac{a^{1/2}}{\pi} \sum_{k\geq 1}\frac{e^{-X_k}}{k^2} \bl\{\sum_{j=0}^{M-1} (-)^j c_j(\fs) X_k^{-j}+O(X_k^{-M})\br\}\]
as $a\to 0+$.
We note that the first double sum and the sum involving the coefficients $B_j$ in (\ref{e46}) vanish when $\nu=\fs$. This results from the fact that the poles at $s=-1, -2, -3, \ldots$ of the integrand in (\ref{e21}) are cancelled by the trivial zeros of $\zeta(2s)$, with the consequence that there are no poles situated in $\Re (s)<0$ when $\nu=\fs$.
This case corresponds to an example of the generalised Euler-Jacobi series mentioned in (\ref{e12}). It can be verified that the above expansion agrees with that given in Theorem 1 of \cite{P16}. 
\vspace{0.6cm}

\begin{center}
{\bf 4. \  Numerical results and concluding remarks}
\end{center}
\setcounter{section}{4}
\setcounter{equation}{0}
\renewcommand{\theequation}{\arabic{section}.\arabic{equation}}
To carry out a numerical investigation of the expansion (\ref{e46}) we define the quantity ${\hat S}_\nu(a)$ by
\[{\hat S}_\nu(a):=S_\nu(a)-{\cal H}(a;\nu)-\cos \pi\nu \sum_{k\geq1} \sum_{n=N_{n-1}}^{N_n-1} \frac{(4n-4\nu)!}{n! (n-\nu)!} \bl(\frac{a}{8\pi^2}\br)^{\!2n-2\nu}\!\zeta(4n\!-\!4\nu\!+\!1,k),\]
where we have used (\ref{e46a}) to express the double sum in alternative form. The computed value of ${\hat S}_\nu(a)$ for a given $a$ and $\nu$ is then compared with the value of the asymptotic sum
\bee\label{e41}
\frac{2^{\nu-\frac{1}{2}}a^\nu}{\pi^{2\nu}}\sum_{k\geq 1} \frac{e^{-X_k}}{k^{1+2\nu}}\bl\{\frac{\cos \pi\nu}{\sqrt{2\pi X_k}} \sum_{j=0}^{M-1}(-)^j B_j X_k^{-j}-\sin \pi\nu \sum_{j=0}^{M-1} (-)^j c_j(\nu) X_k^{-j}\br\}
\ee
for different truncation index $M$. In this comparison we take $k=1$; that is, we consider only the leading subdominant exponential $e^{-X_1}$.

The coefficients $c_j(\nu)$ are obtained from (\ref{e33a}); we recall that the coefficients $B_j=B_j(\alpha,\nu)$
and so must be computed for each value of $a$ (with $\nu$ fixed). As an example in the case $\nu=0$, $a=0.10$,
we have $X_1=10\pi^2$ and, from (\ref{e32}), $N_1=49$, $\alpha=0.69604\dots\ $. The coefficients $B_j$ are computed from (\ref{e45a}) and (\ref{e45}), with the quantity $\gamma_j$ appearing in the coefficients $A_{k,j}$ and $G_{2k,j}$ given by $\gamma_j=-\alpha-j$. An example of these two sets of coefficients is presented in Table 1. 
\begin{table}[th]
\caption{\footnotesize{Values of the coefficients $B_j=B_j(\alpha,\nu)$ and $c_j(\nu)$ for $0\leq j\leq 4$ when $a=0.10$ for two values of the parameter $\nu$.}} \label{t1}
\begin{center}
\begin{tabular}{|l|ll|ll|}
\hline
&&&&\\[-0.25cm]
\mcol{1}{|c|}{} & \mcol{2}{c|}{$\nu=0$} & \mcol{2}{c|}{$\nu=\f{1}{6}$}\\
\mcol{1}{|c|}{} & \mcol{2}{c|}{$N_1=49,\ \alpha=0.6960440109$} & \mcol{2}{c|}{$N_1=50,\ \alpha=-0.8039559891$} \\
\mcol{1}{|c|}{$j$} & \mcol{1}{c}{$B_j$} & \mcol{1}{c|}{$c_j(\nu)$} & \mcol{1}{c}{$B_j$} & \mcol{1}{c|}{$c_j(\nu)$}\\
&&&&\\[-0.3cm]
\hline
&&&&\\[-0.3cm]
0 & $+3.7254213351$ & $1.0000000000$ & $+3.7254213351$ & $1.0000000000$ \\
1 & $-0.4718731128$ & $0.3750000000$ & $+1.1980414492$ & $0.6666666667$ \\
2 & $+1.7780006019$ & $0.4453125000$ & $+3.7772677362$ & $1.0648148148$ \\
3 & $+4.2099300410$ & $0.9228515625$ & $+12.673014554$ & $2.7160493827$ \\
4 & $+13.575274633$ & $2.7781677246$ & $+50.144293786$ & $9.5721882426$ \\
[.1cm]
\hline
\end{tabular}
\end{center}
\end{table}

In Table 2 we show the absolute relative error in the computation of the expansion in (\ref{e41}) as a function of the truncation index $M$. Finally, Table 3 presents the values of ${\hat S}_\nu(a)$ and the asymptotic series (\ref{e41}) with $k=1$ and truncation index $M=5$ for different values of the parameters $a$ and $\nu$. It is seen that there is very good agreement between these values and that this agreement improves as $a$ decreases. The value of the second exponential $e^{-X_2}=e^{-4\pi^2/a}$ in (\ref{e41}) is many orders of magnitude smaller than the first exponential $e^{-X_1}$ and so can safely be ignored at the $k=1$ level.

The situation when only the first exponential $e^{-X_1}$ is taken into account in the expansion (\ref{e41})
is straightforward: the sum is truncated at some suitable point thereby introducing a truncation error. If truncation is optimal, then the resulting error is exponentially more recessive than the exponential factor $e^{-X_1}$. However, when we attempt to include the second exponential corresponding to $k=2$ the situation is not so obvious. It is not apparent, without further investigation, how the error from the optimally truncated leading exponential series compares with the contribution from the second series.
\begin{table}[th]
\caption{\footnotesize{Values of the absolute relative error in the expansion in (\ref{e41}) (with $k=1$) for different truncation index $M$ when $a=0.10$.}} \label{t2}
\begin{center}
\begin{tabular}{|l|ll|}
\hline
&&\\[-0.25cm]
\mcol{1}{|c|}{$M$} & \mcol{1}{c}{$\nu=0$} & \mcol{1}{c|}{$\nu=\f{1}{6}$} \\
&&\\[-0.3cm]
\hline
&&\\[-0.3cm]
1 & $1.329\times 10^{-3}$ & $5.587\times 10^{-3}$ \\
2 & $4.779\times 10^{-5}$ & $8.031\times 10^{-5}$ \\
3 & $1.137\times 10^{-6}$ & $1.962\times 10^{-6}$ \\
4 & $3.671\times 10^{-8}$ & $6.580\times 10^{-8}$ \\
5 & $1.639\times 10^{-9}$ & $2.792\times 10^{-9}$ \\
[.1cm]
\hline
\end{tabular}
\end{center}
\end{table}
\begin{table}[th]
\caption{\footnotesize{Values of ${\hat S}_\nu(a)$ compared with the asymptotic value (\ref{e41}) (with $k=1$ and truncation index $M=5$) for different values of the parameters $a$ and $\nu$. The values of the corresponding optimal truncation index $N_1$ and $\alpha$ are shown.}} \label{t3}
\begin{center}
\begin{tabular}{|lll|ll|}
\hline
&&&&\\[-0.35cm]
\mcol{3}{|c|}{$\nu=0$} & \mcol{2}{c|}{}\\
\mcol{1}{|c}{$a$} & \mcol{1}{c}{$N_1$} & \mcol{1}{c|}{$\alpha$} & \mcol{1}{c}{${\hat S}_\nu(a)$} & \mcol{1}{c|}{Asymptotic}\\
\hline
&&&&\\[-0.3cm]
0.50 & $10$ & $-0.2607911978$ & $+2.9705500754\times 10^{-10}$ & $+2.9705{\bf 7}73493\times 10^{-10}$ \\
0.20 & $25$ & $-0.6519779946$ & $+1.4926163194\times 10^{-23}$ & $+1.492616{\bf 5}957\times 10^{-23}$ \\
0.10 & $49$ & $+0.6960440109$ & $+1.4516136827\times 10^{-44}$ & $+1.45161368{\bf 5}0\times 10^{-44}$ \\
[.1cm]
\hline
&&&&\\[-0.35cm]
\mcol{3}{|c|}{$\nu=\frac{1}{6}$} & \mcol{2}{c|}{}\\
\mcol{1}{|c}{$a$} & \mcol{1}{c}{$N_1$} & \mcol{1}{c|}{$\alpha$} & \mcol{1}{c}{${\hat S}_\nu(a)$} & \mcol{1}{c|}{Asymptotic}\\
\hline
&&&&\\[-0.3cm]
0.50 & $10$ & $+0.2392088022$ & $-3.5340060526\times 10^{-10}$ & $-3.5340{\bf 1}47576\times 10^{-10}$ \\
0.20 & $25$ & $-0.1519779946$ & $-6.0750652224\times 10^{-23}$ & $-6.075065{\bf 6}691\times 10^{-23}$ \\
0.10 & $50$ & $-0.8039559891$ & $-2.3888643465\times 10^{-44}$ & $-2.3888643{\bf 5}31\times 10^{-44}$ \\
[.1cm]\hline
\end{tabular}
\end{center}
\end{table}

\vspace{0.6cm}

\begin{center}
{\bf Appendix A:  Determination of the coefficients $c_j(\nu)$ in the expansion (\ref{e3e})}
\end{center}
\setcounter{section}{1}
\setcounter{equation}{0}
\renewcommand{\theequation}{\Alph{section}.\arabic{equation}}
An algorithm for the determination of the coefficients $c_j(\nu)$ in the inverse factorial expansion (\ref{e3e}) is described in \cite[\S 4]{P16}, \cite[p.~46]{PK}. This procedure is found to work well for numerical values of $\nu$ where up to 100 coefficients have been calculated.

An alternative procedure is described below. For convenience, we replace the variable $s$ in (\ref{e3e}) by $2x$ to obtain
\[\frac{4\g(4x-4\nu)}{x \g(x) \g(x-\nu) \g(2x-3\nu)}=\frac{2^{6x-5\nu+\frac{3}{2}}}{2\pi}\bl\{c_0(\nu)-\frac{c_1(\nu)}{2x-3\nu-1}+\frac{c_2(\nu)}{(2x-3\nu-1)(2x-3\nu-2)}+\cdots\br\}.\]
Use of the multiplication formula for the gamma function \cite[(5.5.6)]{DLMF} shows that the left-hand side of the above expression reduces to
\[\frac{2^{6x-5\nu+\frac{3}{2}}}{2\pi x}\,\frac{\g(x-\nu+\f{1}{4}) \g(x-\nu+\fs) \g(x-\nu+\f{3}{4})}{\g(x-\f{3}{2}\nu) \g(x) \g(x-\f{3}{2}\nu+\fs)}.\]
From the asymptotic expansion of the ratio of two gamma functions \cite[(5.11.13)]{DLMF} we have for $x\to+\infty$
\begin{eqnarray*}
{\cal S}_1&=&\frac{\g(x-\nu+\f{1}{4})}{\g(x-\f{3}{2}\nu)}=x^{\frac{1}{2}\nu+\frac{1}{4}}\bl\{1-\frac{3+16\nu+20\nu^2}{32x}+O(x^{-2})\br\} \\ 
{\cal S}_2&=&\frac{\g(x-\nu+\fs)}{\g(x)}=x^{\frac{1}{2}-\nu}\bl\{1+\frac{4\nu^2-1}{8x}+O(x^{-2})\br\} \\ 
{\cal S}_3&=&\frac{\g(x-\nu+\f{3}{4})}{\g(x-\f{3}{2}\nu+\fs)}=x^{\frac{1}{2}\nu+\frac{1}{4}}\bl\{1+\frac{1-8\nu-20\nu^2}{32x}+O(x^{-2})\br\}. 
\end{eqnarray*}

Then
\[\frac{1}{x}\,{\cal S}_1 {\cal S}_2 {\cal S}_3=1-\frac{3(1+2\nu)^2}{8 (2x)}+\frac{9+8\nu-72\nu^2-32\nu^3+144\nu^4}{128 (2x)^2}+\cdots\]
\[\hspace{0.3cm}=c_0(\nu)-\frac{c_1(\nu)}{2x}+\frac{c_2(\nu)-(1+3\nu)c_1(\nu)}{(2x)^2}+\cdots\ ,\]
from which it is possible to compute $c_1(\nu)$ and $c_2(\nu)$. Continuation of this process with the aid of {\it Mathematica} then yields the coefficients
\[c_0(\nu)=1,\quad c_1(\nu)=\frac{3}{8}(1+2\nu)^2,\]
\[c_2(\nu)=\frac{1}{128}(57+344\nu+696\nu^2+544\nu^3+144\nu^4),\]
\[c_3(\nu)=\frac{1}{1024}(945+7068\nu+19196\nu^2+24288\nu^3+15472\nu^4+4800\nu^5+576\nu^6).\] 
\[c_4(\nu)=\frac{1}{32768}(91035+780624\nu+2531344\nu^2+4113984\nu^3+3746848\nu^4+1994496\nu^5\]
\[\hspace{6cm}+615040\nu^6+101952\nu^7+6912\nu^8).\]
\vspace{0.6cm}

\begin{center}
{\bf Appendix B: A rearrangement of a double sum}
\end{center}
\setcounter{section}{2}
\setcounter{equation}{0}
\renewcommand{\theequation}{\Alph{section}.\arabic{equation}}
In order to detect the exponential terms $e^{-X_k}$ in computations it is necessary to remove all larger terms in the first asymptotic series on the right-hand side of (\ref{e46}). This can be achieved by a straightforward regrouping of the terms in the absolutely convergent double sum
\[\sum_{k\geq 1}\frac{1}{k}\sum_{n=m+1}^{N_k-1} \frac{(4n-4\nu)!}{n! (n-\nu)!}\,\bl(\frac{a}{8\pi^2 k^2}\br)^{\!2n-2\nu}\equiv \sum_{k\geq 1} \sum_{n=m+1}^{N_k-1}{\cal F}_{n,k},\]
where
\[{\cal F}_{n,k}=\frac{A_n}{k^{4n+1-4\nu}},\qquad A_n=\frac{(4n-4\nu)!}{n! (n-\nu)!}\,\bl(\frac{a}{8\pi^2 }\br)^{\!2n-2\nu}.\]
Thus we have (with $N_0=m+1$)
\[\sum_{k\geq 1} \sum_{n=m+1}^{N_k-1}{\cal F}_{n,k}=\sum_{n=N_0}^{N_1-1}\sum_{k\geq 1} {\cal F}_{n,k}+\sum_{r=N_1}^{N_2-1} \sum_{k\geq 2} {\cal F}_{n,k}+
\sum_{r=N_2}^{N_3-1} \sum_{k\geq 3} {\cal F}_{n,k}+\ldots\]
\[=\sum_{p=1}^\infty \sum_{n=N_{n-1}}^{N_n-1}\,{\cal A}_n \,\zeta(4n+1-4\nu,p),\]
where $\zeta(x,p)$ is the Hurwitz zeta function.

\end{document}